\documentclass[11pt]{article}
\usepackage{newlfont,amsfonts,amssymb, oldgerm,amsmath,amsthm,amsgen,amscd,yhmath}
\usepackage{epic}
\usepackage{eepic}
\usepackage{rotating}
\usepackage[dvips]{epsfig}
\begin{document}

\newcommand{\Q}{\ensuremath{\mathbb{H}}}
\newcommand{\N}{\ensuremath{\mathbb{N}}}
\newcommand{\Z}{\ensuremath{\mathbb{Z}}}
\newcommand{\C}{\ensuremath{\mathbb{C}}}
\newcommand{\K}{\ensuremath{\mathbb{K}}}
\renewcommand{\O}{\ensuremath{\mathcal{O}}}
\newcommand{\R}{\ensuremath{\mathbb{R}}}

\newcommand{\bcase}{\begin{case}}
\newcommand{\ecase}{\end{case}}
\newcommand{\setcase}{\setcounter{case}{0}}
\newcommand{\bclaim}{\begin{claim}}
\newcommand{\eclaim}{\end{claim}}
\newcommand{\setclaim}{\setcounter{claim}{0}}
\newcommand{\bstep}{\begin{step}}
\newcommand{\estep}{\end{step}}
\newcommand{\setstep}{\setcounter{step}{0}}
\newcommand{\bhlem}{\begin{hlem}}
\newcommand{\ehlem}{\end{hlem}}
\newcommand{\sethlem}{\setcounter{hlem}{0}}

\newcommand{\bleer}{\begin{leer}}
\newcommand{\eleer}{\end{leer}}
\newcommand{\bde}{\begin{de}}
\newcommand{\ede}{\end{de}}
\newcommand{\ul}{\underline}
\newcommand{\ol}{\overline}
\newcommand{\tbf}{\textbf}
\newcommand{\mc}{\mathcal}
\newcommand{\mb}{\mathbb}
\newcommand{\mf}{\mathfrak}
\newcommand{\bs}{\begin{satz}}
\newcommand{\es}{\end{satz}}
\newcommand{\btheo}{\begin{theo}}
\newcommand{\etheo}{\end{theo}}
\newcommand{\bfolg}{\begin{folg}}
\newcommand{\efolg}{\end{folg}}
\newcommand{\blem}{\begin{lem}}
\newcommand{\elem}{\end{lem}}
\newcommand{\bnote}{\begin{note}}
\newcommand{\enote}{\end{note}}
\newcommand{\bprf}{\begin{proof}}
\newcommand{\eprf}{\end{proof}}
\newcommand{\bd}{\begin{displaymath}}
\newcommand{\ed}{\end{displaymath}}
\newcommand{\be}{\begin{eqnarray*}}
\newcommand{\ee}{\end{eqnarray*}}
\newcommand{\eeqa}{\end{eqnarray}}
\newcommand{\beqa}{\begin{eqnarray}}
\newcommand{\bi}{\begin{itemize}}
\newcommand{\ei}{\end{itemize}}
\newcommand{\bnum}{\begin{enumerate}}
\newcommand{\enum}{\end{enumerate}}
\newcommand{\la}{\langle}
\newcommand{\ra}{\rangle}
\newcommand{\eps}{\epsilon}
\newcommand{\ve}{\varepsilon}
\newcommand{\vp}{\varphi}
\newcommand{\lra}{\longrightarrow}
\newcommand{\Lra}{\Leftrightarrow}
\newcommand{\Ra}{\Rightarrow}
\newcommand{\sub}{\subset}
\newcommand{\ems}{\emptyset}
\newcommand{\sms}{\setminus}
\newcommand{\ints}{\int\limits}
\newcommand{\sums}{\sum\limits}
\newcommand{\lims}{\lim\limits}
\newcommand{\bcup}{\bigcup\limits}
\newcommand{\bcap}{\bigcap\limits}
\newcommand{\beq}{\begin{equation}}
\newcommand{\eeq}{\end{equation}}
\newcommand{\einhalb}{\frac{1}{2}}
\newcommand{\rr}{\mathbb{R}}
\newcommand{\rn}{\mathbb{R}^n}
\newcommand{\ccc}{\mathbb{C}}
\newcommand{\cn}{\mathbb{C}^n}
\newcommand{\M}{{\cal M}}
\newcommand{\drehgleich}{\mbox{\begin{rotate}{90}$=$  \end{rotate}}}
\newcommand{\turngleich}{\mbox{\begin{turn}{90}$=$  \end{turn}}}
\newcommand{\turnsimeq}{\mbox{\begin{turn}{270}$\simeq$  \end{turn}}}
\newcommand{\vf}{\varphi}
\newcommand{\earr}{\end{array}\]}
\newcommand{\barr}{\[\begin{array}}
\newcommand{\bvec}{\left(\begin{array}{c}}
\newcommand{\evec}{\end{array}\right)}
\newcommand{\sumk}{\sum_{k=1}^n}
\newcommand{\sumi}{\sum_{i=1}^n}
\newcommand{\suml}{\sum_{l=1}^n}
\newcommand{\sumj}{\sum_{j=1}^n}
\newcommand{\sumij}{\sum_{i,j=1}^n}
\newcommand{\sumkl}{\sum_{k,l=1}^n}
\newcommand{\suminf}{\sum_{k=0}^\infty}
\newcommand{\inv}{\frac{1}}
\newcommand{\wzbw}{\hfill $\Box$\\[0.2cm]}
\newcommand{\lag}{\mathfrak{g}}
\newcommand{\lan}{\mathfrak{n}}
\newcommand{\lah}{\mathfrak{h}}
\newcommand{\lak}{\mathfrak{k}}
\newcommand{\lam}{\mathfrak{m}}
\newcommand{\laz}{\mathfrak{z}}
\newcommand{\+}{\oplus}
\newcommand{\x}{\times}
\newcommand{\lx}{\ltimes}
\newcommand{\rrn}{\mathbb{R}^n}
\newcommand{\laso}{\mathfrak{so}}
\newcommand{\lason}{\mathfrak{so}(n)}
\newcommand{\lagl}{\mathfrak{gl}}
\newcommand{\lasl}{\mathfrak{sl}}
\newcommand{\lasp}{\mathfrak{sp}}
\newcommand{\lasu}{\mathfrak{su}}
\newcommand{\w}{\omega}
\newcommand{\pmh}{{\cal P}(M,h)}
\newcommand{\s}{\sigma}
\newcommand{\deri}{\frac{\partial}}
\newcommand{\ddx}{\frac{\partial}{\partial x}}
\newcommand{\ddz}{\frac{\partial}{\partial z}}
\newcommand{\ddi}{\frac{\partial}{\partial y_i}}
\newcommand{\ddj}{\frac{\partial}{\partial y_j}}
\newcommand{\ddk}{\frac{\partial}{\partial y_k}}
\newcommand{\ddp}{\frac{\partial}{\partial p_i}}
\newcommand{\ddq}{\frac{\partial}{\partial q_i}}
\newcommand{\ddl}{\frac{\partial}{\partial y_l}}
\newcommand{\xz}{^{(x,z)}}
\newcommand{\mh}{(M,h)}
\newcommand{\wxz}{W_{(x,z)}}
\newcommand{\qmh}{{\cal Q}(M,h)}
\newcommand{\bbem}{\begin{bem}}
\newcommand{\ebem}{\end{bem}}
\newcommand{\bbez}{\begin{bez}}
\newcommand{\ebez}{\end{bez}}
\newcommand{\bbsp}{\begin{bsp}}
\newcommand{\ebsp}{\end{bsp}}
\newcommand{\pr}{pr_{\lason}}
\newcommand{\huts}{\hat{\s}}
\newcommand{\whut}{\w^{\huts}}
\newcommand{\bhg}{{\cal B}_H(\lag)}
\newcommand{\aaa}{\alpha}
\newcommand{\bb}{\beta}
\newcommand{\laa}{\mf{a}}
\newcommand{\Ll}{\lambda}
\newcommand{\LL}{\Lambda}
\newcommand{\D}{\Delta}
\newcommand{\ß}{\beta}
\newcommand{\ä}{\alpha}
\newcommand{\W}{\Omega}
\newcommand{\esel}{\ensuremath{\mathfrak{sl}(2,\ccc)}}
\newcommand{\kg}{{\cal K}(\lag)}
\newcommand{\bg}{{\cal B}_h(\lag)}
\newcommand{\kk}{  \mathbb{K}}
\newcommand{\xy}{[x,y]}
\newcommand{\perdef}{$\stackrel{\text{\tiny def}}{\iff}$}
\newcommand{\eqdef}{\stackrel{\text{\tiny def}}{=}}
\newcommand{\lai}{\mf{i}}
\newcommand{\lar}{\mf{r}}
\newcommand{\Dim}{\mathsf{dim\ }}
\newcommand{\im}{\mathsf{im\ }}
\newcommand{\Ker}{\mathsf{ker\ }}
\newcommand{\trace}{\mathsf{trace\ }}
\newcommand{\grad}{\mathsf{grad}}
\newcommand{\lecturecount}{\begin{center}{\sf  (Lecture \Roman{lecturenr})}\end{center}\addtocontents{toc}{{\sf  (Lecture \Roman{lecturenr})}} \refstepcounter{lecturenr}
}
\newcommand{\T}{{\cal T}}
\newcommand{\cur}{{\cal R}} 
\newcommand{\pd}{{\cal P}}  

\theoremstyle{definition}
\newtheorem{de}{Definition}
\newtheorem{bem}{Remark}
\newtheorem{bez}{Notation}
\newtheorem{bsp}{Example}
\theoremstyle{plain}
\newtheorem{lem}{Lemma}
\newtheorem{satz}{Proposition}
\newtheorem{folg}{Corollary}
\newtheorem{theo}{Theorem}

\bibliographystyle{alpha}


\title{Screen bundles of Lorentzian manifolds and some generalisations of pp-waves}

\author{Thomas Leistner}
\date{ }
\maketitle

\begin{abstract}
A pp-wave is a Lorentzian manifold with a parallel light-like vector field satisfying a certain curvature condition. We introduce  generalisations of pp-waves, on one hand by allowing the vector field to be recurrent and on the other hand by weakening the curvature condition. These generalisations are related to the screen holonomy of the Lorentzian manifold. While pp-waves have a trivial screen holonomy there are no restrictions on the screen holonomy  of the manifolds with the weaker curvature condition. 
\\[.2cm]
\end{abstract}


\section*{Introduction}

Regarding holonomy theory or the existence of parallel spinors, undoubtfully
the most interesting Lorentzian manifolds are those with indecomposable,
but non-irreducible holonomy representation.
They admit a recurrent light-like vector field and their holonomy algebra is
contained in the parabolic algebra $(\rr \+ \lason)\ltimes \rrn$, assumed that the dimension of the manifold is $n+2$. The main ingredient of this holonomy algebra is its $\lason$-projection, which is called {\em screen holonomy}.
In previous papers  
 \cite{leistner02}, \cite{leistner03}, \cite{leistner03b} 
we adressed ourselves to the classification of the screen holonomy and obtained
the result that it has to be a Riemannian holonomy.

On the other hand it can be shown that any Riemannian holonomy group can be realised as screen holonomy of an indecomposable, non-irreducible Lorentzian manifold by a rather simple method: for a Riemannian manifold $(N,g)$ and $f\in C^\infty(N)$ the manifold $\rr^2\times N$ with Lorentzian metric $2dxdz+fdz^2 +g$ is non-irreducible, indecomposable for $f$ sufficiently generic  and, above all, its screen holonomy is equal to the Riemannian holonomy of $(N,g)$.

In this note we want to consider Lorentzian manifolds which are in some sense complementary to the ones obtained by this procedure and which can be understood as certain generalisations of pp-waves. pp-waves are defined by the existence of a light-like parallel vector field and a certain curvature condition. Or aim is to generalise pp-waves in two directions: on one hand we will only require the existence of a recurrent vector field instead of a parallel one (see Section \ref{secpr}), and on the other hand, more importantly we will relax the curvature condition (see Section \ref{secrr}). These generalisations are related to the screen holonomy in the following sense. pp-waves have trivial screen holonomy, i.e. their screen bundle, which we will  introduce in Section \ref{secscreen} is flat. This remains true if we drop the assumption that the vector field is parallel, but it is no longer true if we weaken the curvature condition. Instead, one can prove that the screen bundles restricted to the light-like hypersurfaces defined by the recurrent vector field are flat.

These generalisations can also be understood in terms of the ingredients of the local form of a Lorentzian metric $h$ with recurrent vector field which are a function $f$, a $1$-form $\phi$ and a family of Riemannian metrics $g_z$ because  $h$  can be written as 
$h=2dxdz + fdz^2 + \phi dz + g_z$. For a pp-wave it is $\phi=0$, $g_z$ flat and $\ddx(f)=0$. If we no longer require a parallel vector field only the conditions $\phi=0$ and $g_z$ flat remain. Finally, weakening the curvature conditions is equivalent to dropping the assumption $\phi=0$, i.e. only requiring $g$ to be flat. As mentioned this is complementary to the construction method above where $\phi=0$ is obtained.

Although the curvature conditions to these generalised pp-waves are only slightly weaker the consequences for the screen holonomy are dramatical in the following sense. While for pp-waves the screen holonomy has to be trivial, any possible screen holonomy, that is any Riemannian holonomy, can be obtained for the generalisations of pp-waves. This can be deduced from a recent result of Galaev in \cite{galaev05} and is explained in the last section.
\section{Lorentzian manifolds with recurrent  light-like vector field}
\label{seceins}
A vector field $X$ is called recurrent if $\nabla X = \Theta \otimes X$ where 
$\Theta $ is a one-form on $M$. If the lenght of a recurrent vector field is non-zero, it can be rescaled to a parallel one. This is not true in general if
the recurrent vector field is lightlike. 

If a Lorentzian manifold $(M,h)$ carries a recurrent light-like vector field $X$  the holonomy group of $(M,h)$ in $p\in M$ admits the one-dimensional light-like invariant subspace $\rr\cdot X_p$, hence it does not act irreducible. The orthogonal complement of this subspace $X_p^\bot$ is $n+1$-dimensional, holonomy invariant as well and contains $\rr\cdot X_p$. Hence $X$ yields two parallel distributions, a one-dimensional, totally isotropic distribtution  $\Xi$ with $X\in \Gamma (\Xi)$, and its $n+1$ dimensional orthogonal complement   $\Xi^\bot=\{U\in TM\mid h(U,X)=0\}$ containing $\Xi$. Both foliate the manifold into light-like lines $\cal X$, which are the flow of $X$, and light-like hypersurfaces $\cal X^\bot$. Using this foliation the following coordinate description was proven (see \cite{walker49}, \cite{brinkmann25} and \cite{schimming74}).

\bs
Let $(M,h)$ be a Lorentzian manifold of dimension $n+2>2$  with recurrent vector field $X$.
\bnum
\item This is equivalent to the existence of 
 coordinates $(U,\varphi=(x, (y_i)_{i^=1}^n, z)) $ in which the metric $h$ has the following
local shape
\begin{equation}
\label{walker}h = 2\ dx dz + \sum_{i = 1}^n u_i dy_i dz  + f dz^2 +  \sum_{i,j =
1}^n g_{ij} dy_i\ dy_j  \end{equation}
with $ \frac{\partial g_{ij}}{\partial x}= \frac{\partial u_i}{\partial x}=0$, $f\in
C^\infty(M)$. To these coordinates we refer as {\em Walker coordinates}.
\item
$X$ is parallel if and only if $f$ does not depend on $x$.  To these coordinates we refer as
 {\em Brinkmann coordinates}.
\item 
If $X$ is parallel the coordinates can be chosen such that $u_i=0$ and end even that $f=0$. 
 To these coordinates we refer as {\em Schimming coordinates}.
\enum
\es

A Lorentzian manifold with lightlike parallel vector field is called {\em
Brinkmann-wave}, after \cite{brinkmann25}. For further coordinate descriptions see \cite{boubel00} or \cite{thesis}.

Returning to the holonomy group of a Lorentzian manifold with recurrent light-like vector field we want to mention some of its algebraic properties.
The holonomy algebra $\mf{h}$ of a $n+2$-dimensional  Lorentzian manifold with recurrent vector field is contained in 
the parabolic algebra $\mf{p}=(\rr\+\lason )\ltimes \rrn$ which is given in an appropriate basis as
\[\mf{p}
=\left\{
\left.
\left(\begin{array}{ccc}a & v^t & 0 \\0 & A & v \\0 & 0 & -a\end{array}\right)\right|
a\in \rr, v\in \rrn, A\in \lason\right\}.\]
Its projection onto $\rrn$ is surjective if and only if the holonomy group acts  indecomposably. The recurrent vector field is parallel if and only if the holonomy is
contained in $\lason\ltimes \rrn$. 
There are four different algebraic types of holonomy algebras (see \cite{bb-ike93}),
two of them {\em uncoupled}, i.e. $\mf{h}=\lag\ltimes \rrn$ and $\mf{h}=(\rr\+\lag)\ltimes \rrn$, and two with a coupling
between the center of the $\lason$-projection and the $\rr$-- resp. the $\rrn$--part. 
 Further algebraic properties can be proved easily.
\blem
Let  $\mf{h}$ be an  indecomposable  subalgebra  of the parabolic algebra. Then:
\bnum
\item
$\mf{h}$ is solvable if and only if it is 2-step solvable (i.e. $\mf{h}^{(1)}\not= 0$  and $\mf{h}^{(2)}= 0$) or Abelian.
\item
$\mf{h}$ is Abelian if and only if  $\mf{h}=\rrn$.
\item If 
$\mf{h}$ is the holonomy algebra of an indecomposable, non-irreducible Lorentzian manifold, then
$\mf{h}$ is 2-step solvable if and only if $\mf{h}=\rr\ltimes\rrn$ or the screen holonomy algebra equals to a direct sum
of  copies of $\laso (2)$.
\enum 
\elem
\bprf
The first point is obvious from the commutator relations in the parabolic algebra
\[ [(a,A,x),(b,B,y)] \ =\ \left(0,[A,B], (A+a Id)y - (B+bId)x\right).\]
Set $\lag:=pr_{\lason}\mf{h}$. If $\lah$ is solvable, then $\lag$ has to be solvable. But, as a subalgebra 
of $\lason$, $\lag$ is reductive, i.e. it is solvable if and only if it is Abelian. Hence
$\lah^{(1)}\subset \lag^{(1)}\ltimes \rrn = \rrn$ and therefore $\lah^{(2)}=0$.  From the commutator relation one sees 
that $\lah^{(1)}=0$ only if $\lah=\rrn$. 
The remaining decomposition of the $\lason$--part $\lag$ under the assumption that $\lah$ is a holonomy algebra
follows from a Borel-Lichnerowicz decomposition theorem proved in \cite{bb-ike93}.
\eprf

\section{The screen bundle associated to a recurrent vector field}
\label{secscreen}
In this section we will describe the $SO(n)$--projection of an indecomposable, non-irreducible
holonomy group of a $n+2$--dimensional, simply connected Lorentzian manifold
as a holonomy group of a metric connection in a vector bundle, the so called
{\em screen bundle}. The results of this section were obtained in \cite{thesis}.

We consider  the distributions $\Xi$ and $\Xi^\bot$ on $M$ introduced in Section \ref{seceins}, which are parallel, i. e.
$\nabla_U$ leaves $\Gamma(\Xi)$ and $\Gamma(\Xi^\bot)$ invariant for all $U\in TM$.
The factor spaces $\Xi^\bot_p/\Xi_p$ in every point $p\in M$  define a
vector bundle over M,
\be
{\cal S}&:=&\dot{\bigcup_{p\in M}} \Xi_p^\bot/\Xi_p,
\ee
which is  called {\em screen bundle}.
The
metric $h$ on $M$ defines a scalar product on ${\cal S}$, which we denote by $\hat{h}$, via
\[ \hat{h}\left( [X], [Y] \right) := h (X,Y).\]
With respect to this scalar product the bundle $ {\cal O}({\cal S})$ is defined as
the set of  orthonormal frames of ${\cal S}$ over $M$. This is a $ O
(n)$--principal fibre bundle. ${\cal O}({\cal S})$ has fibres
\[
{\cal O}_p({\cal S})= \left\{ \left( [E_1], \ldots , [E_n] \right) \left|
\begin{array}{l}
 (X, E_1, \ldots , E_n )
\mbox{ a basis of }\Xi_p^\bot \mbox{ for } X\in \Xi_p\\
\mbox{with }h(E_i, E_j) = \delta_{ij}
\end{array}\right. \right\}.\]
Then we can describe ${\cal S}$ as vector bundle associated to the bundle ${\cal O}({\cal S})$:
\begin{eqnarray*}
 {\cal O}({\cal S}) \times_{O (n)} \mathbb{R}^n&\simeq& {\cal S}
\\
\left[ \left([E_1], \ldots , [E_n] \right), ( x_1, \ldots x_n) \right]
&\mapsto&
\left[ \sum_{i = 1}^n x_i E_i \right]
\end{eqnarray*}

We now consider subbundle  ${\cal P}(M,h)$ of the frame bundle with fibres
\beq{\cal P}_p (M,h):=\left\{(X, E_1 , \ldots E_n,Z) \left|
\begin{array}{l} X\in \Xi_p, E_i\in \Xi_p^\bot,  h(E_i,E_j)= \delta_{ij},\\
h(Z,Z)=h(Z,E_i)=0, h(X,Z)=1
\end{array}\right.
	\right\} \label{bundlep}
\eeq and structure group $P=(\rr^*\times O(n))\ltimes \rrn$. We define a surjective bundle homomorphism
\[\begin{array}{rcccc}
f&:& {\cal P}(M,h) & \rightarrow & {\cal O}({\cal S})\\
&& \left(X, E_1, \ldots , E_n, Z \right) & \mapsto & \left( [E_1], \ldots , [E_n]
\right).
\end{array}\]
Then $f$ defines a reduction of the projection $pr_{O(n)}: P=(\rr^*\times O(n))\ltimes \rrn\rightarrow O(n)$.
\blem
$f:{\cal P}(M,h)  \rightarrow  {\cal O}({\cal S})$ is a $pr_{O (n)}$--reduction.
\elem
\bprf
We have to verify that the following
diagram commutes
\barr{ccccl}
P\times {\cal P}(M,h)&\longrightarrow &{\cal P}(M,h)&&\\
&&&\searrow&\\
\makebox[0cm][r]{{\scriptsize  $pr_{SO(n)}\times f$}}\downarrow&\circlearrowright&
\makebox[0cm][r]{{\scriptsize $f$}}\downarrow\makebox[0cm][l]{ $\;\;\circlearrowright$}&&M.\\
&&&\nearrow&\\
O(n)\times {\cal O}({\cal S})&\longrightarrow &{\cal O}({\cal S})&&
\earr
The action of the components of $P$ on $ {\cal P}(M,h)$ is as follows:
\beq\label{rsternact}
(X, E_1, \ldots , E_n, Z)\cdot (a,Id,0)=(aX, E_1, \ldots , E_n, a^{-1}Z)
\eeq
and
\beq\label{vectoract}
\begin{array}{l}
(X, E_1, \ldots , E_n, Z)\cdot (1,Id,v)=
\\
(X, v_1 X+E_1 , \ldots , v_nX+E_n, - \frac{1}{2} v^t v\  X - \sum\limits_{k=1}^n v_k E_k +Z ).
\end{array}
\eeq
Since $P$ is a semi-direct product this implies that
\[ f\left( (X, E_1, \ldots , E_n, Z)\cdot (a,A,v)\right)=([E_1], \ldots , [E_n] )\cdot A.\]
But this makes the   the above diagram commutative.
\eprf

Since $\Xi$ is parallel the Levi-Civita connection defines also a covariant
derivative $\nabla^{{\cal S}}$ on ${\cal S}$ by
\[ \nabla^{{\cal S}}_X [Y] := \left[ \nabla_X Y \right]. \]
This covariant derivative is metric with respect to $\hat{h}$ since the Levi-Civita
connection is metric.
It defines a  connection form $ \theta $ on $ {\cal O}({\cal S})$ which is given for a local
section $ \hat{\sigma}= \left( [\sigma_1], \ldots ,[ \sigma_n] \right)\in\Gamma({\cal O}({\cal S}))$ by the formula
\[ \nabla^{{\cal S}}_U [V] = \nabla^{{\cal S}}_U \left[ ( \hat{\sigma} , \nu ) \right] = \left[ \hat{\sigma} ,
d\nu(V) + \theta^{ \hat{\sigma}} (U) \cdot  \nu \right] \]
for $\nu=(\nu_1,\ldots,\nu_n)$ and $[V]= \sum_{i = 1}^n \nu_i [\sigma_i] $ locally, where $\theta^{\hat{\sigma}}$
is the local connection form of $\theta$.
We get the following result.

\bs
Let $(M,h)$ be an indecomposable Lorentzian manifold of dimension $n+2$ and with parallel isotropic distribution
$\Xi$. Let $\w$ denote the connection form of the Levi-Civita connection $\nabla$.
Then $\w$ is a $pr_{O(n)}$-reduction of the connection $\theta$ of ${\cal O}({\cal S})$.
\es

\bprf
We consider the diagram
\beq\label{reduction}
\begin{CD}
T{\cal P}(M,h)@>{df}>>T{\cal O}({\cal S})\\
@V{\w}VV     @VV{\theta}V\\
\mf{p}@>>{dpr_{O(n)}=pr_{\lason}}>\lason
\end{CD}\eeq
and have to show that $(df)_s$ sends the kernel of
$\w_s$ to the kernel of $\theta_{f(s)}$ for $s\in \pmh$.

Now every element in the kernel of $\w_s$ is equal to
$(d\sigma)_p (U)$ for $p\in M$, $U\in T_pM$ and a certain local section
$\sigma\in \Gamma(\pmh)$ with $\sigma(p)=s$.
Now it is
\[0=\w_{\sigma(p)}((d\sigma)_p(U))=(\sigma^*\w)_p(U)=\w^\sigma_p(U).\]
For the local connection  form  $\w^\sigma$ of the Levi-Civita connection one calculates as follows:
for $\sigma=(\xi, \sigma_1, \ldots , \sigma_n , \zeta )\in \Gamma(\pmh)$ and
${E}_{rt}$ the standard basis of $\mf{gl}(n,\rr)$ it is
\barr{rcll}
0\ =\ \w^\sigma(U)&= &
h(\nabla_U\xi,\zeta)\left( E_{00}-E_{n+1 n+1}\right)&\mbox{(the $\rr$--part)}\\
&&{}+\sums_{k=1}^n h(\nabla_U\sigma_k,\zeta)\left( E_{0k}-E_{k n+1}\right)&\mbox{(the $\rrn$--part})\\
&&{}+\sums_{1\le k< l\le n} h(\nabla_U\sigma_k,\sigma_l)\left( E_{kl}-E_{lk}\right)&\mbox{(the $\lason$--part)}.
\earr
We have to consider $(df)_{\sigma(p)}(d\sigma)_p(U)=d(f\circ \sigma)_p(U)$. If now $\sigma\in \Gamma(\pmh)$ as above,
then is $f\circ \sigma=([\sigma_1], \ldots ,[\sigma_n] )\in \Gamma({\cal O}({\cal S}))$. Finally it is
\be
\theta_{f\circ \sigma(p)}(d(f\circ\sigma)_p(U))&=&\theta^{f\circ \sigma}(U)\\
&=&\sums_{1\le k< l\le n} \hat{h}(\nabla^{\cal S}_U[\sigma_k],[\sigma_l])\left( E_{kl}-E_{lk}\right)\\
&=&\sums_{1\le k< l\le n} h(\nabla_U\sigma_k,\sigma_l)\left( E_{kl}-E_{lk}\right)\\
&=&0
\ee
because of the equation above. I.e. $d(f\circ \sigma)_p(U)$ is in the kernel of the local
connection $\theta^{f\circ \sigma}$. Hence it is in the kernel of $\theta$.
\eprf

\bfolg\label{vbhol}
The diagram (\ref{reduction}) commutes and the curvatures $\Theta$ of $\theta$ and $\Omega$ of $\w $ satisfy
\[f^*\Theta = pr_{\lason}\circ \Omega.\]
This implies the following for the holonomy algebras:
\[ \mf{hol}_p({\cal S},\nabla^{\cal S})=pr_{\lason}(\mf{hol}_p(M,h)).\]
\efolg

\bprf This follows from the proposition by the general theory of reductions of connections.
Since $f$ and $df$ are surjective one gets by the Ambrose-Singer holonomy holonomy theorem that
$\mf{hol}_{f(s)}(\theta)=pr_{\lason}(\mf{hol}_{s}(\w))$.
\eprf

In \cite{leistner02}, \cite{leistner03} and \cite{leistner03b} we have shown that the screen holonomy $\lag:=pr_{\lason}(\mf{hol}(M,h))$ has to be a Riemannian holonomy algebra.
Furthermore, the description of 
$\lag$ as holonomy of the screen bundle
can be used to interpret the geometric information which is algebraically encoded in $\lag$
as geometric structure on the
screen bundle ${\cal S}$. For example, if there is a complex structure on ${\cal S}$ which is
compatible with the metric $\hat{h}$ and parallel to the covariant derivative $\nabla^{\cal S}$ then the flag $\Xi\subset \Xi^\bot\subset TM$ is called
{\em  K\"ahler flag}. The existence of such a K\"ahler flag is equivalent to
$\lag\subset \mf{u}(n)$. For $\lag\subset\mf{su}(n)$ one calls such a flag {\em special K\"ahler flag}.
For details see
\cite{baumsurvey} and \cite{kathhabil}. This can be done analogously for any other geometric structure
on ${\cal S}$, resp. algebraic structure on $\lag$.

\section{Lorentzian manifolds with trivial screen holonomy} 
\label{secpr}
In this section we want to recall results about pp-waves which lead to a further generalisation of pp-waves in the next section.
But first we recall the definition of a pp wave.
A Brinkmann-wave is called {\em $pp$-wave} if its curvature tensor ${\cal R}$
satisfies the
trace condition
$tr_{(3,5)(4,6)} ( {\cal R} \otimes {\cal R} ) =0$.
Schimming  proved the following coordinate description and equivalences in  \cite{schimming74}. 
\blem
A Lorentzian manifold $(M,h)$ of dimension $n+2>2$ is a pp-wave if and only if 
there exist local coordinates
 $(U,\varphi=(x, (y_i)_{i=1}^n, z)) $ in which the metric $h$ has the
form
\begin{equation}
\label{ppform}h = 2\ dx dz   + f dz^2 +  \sum_{i =
1}^n  dy_i^2 \ \mbox{, with $ \frac{\partial f}{\partial x}= 0$.} 
\end{equation}
\elem
\blem\label{eqs}
A Brinkmann wave $(M,h)$ with parallel lightlike vector field $X$ is a pp-wave if and only if
one of the following conditions --- in which $\xi$ denotes the 1-form $h(X,.)$ ---
is satisfied:
\begin{enumerate}
\item \label{eq1}
$ \Lambda_{(1,2,3)} \left(\xi \otimes {\cal R}\right) =0$
\item\label{eq2}
There is a symmetric  $(2,0)$-tensor $\varrho$ with $\varrho(X,.)=0$, such that \\
$ {\cal R} =  \ \Lambda_{(1,2)(3,4)}\left( \xi\otimes \varrho \otimes \xi\right)$.
\item\label{eq3}
There is a function $\vf$, such that
$ tr_{(1,5)(4,8)} ( {\cal R} \otimes {\cal R} ) = \vf\ \xi \otimes\xi  \otimes 
\xi \otimes \xi$.
\end{enumerate}
\elem

In \cite{leistner05} we gave another equivalence for the definition which seems to be simpler than any of
the trace conditions and which makes the generalisation in the next section possible. 

\bs
A Brinkmann-wave $(M,h)$ with parallel lightlike vector field $X$ and induced parallel distributions $\Xi$ and $\Xi^\bot$
 is a pp-wave if and only if
its curvature tensor satisfies:
\begin{equation}
\label{ppeinfach1}
\cur (U,V): \Xi^\bot \longrightarrow \Xi \mbox{ for all }U,V\in TM,
\end{equation}
or equivalently 
\begin{equation}
\label{ppeinfach2}
\cur (Y_1,Y_2)=0 \mbox{ for all } Y_1,Y_2\in \Xi^\bot .
\end{equation}
\es

From this description one obtains easily that a pp-wave is Ricci-isotropic and has vanishing scalar curvature. But it also enables us to introduce a class of non-irreducible Lorentzian manifold by supposing (\ref{ppeinfach1}) but only the existence of a recurrent vector field.
Assuming that the abbreviation `pp' stands for `plane fronted with parallel rays' we call them {\em
pr-waves}: `plane fronted with recurrent rays'.

\begin{de}\label{prdef}
We call a Lorentzian manifold $(M,h)$ {\em 
pr-wave} if it admits a recurrent vector field $X$ and   its curvature tensor ${\cal R}$ obeys
\begin{equation}\label{preq}
\cur (U,V): \Xi^\bot \longrightarrow\Xi \mbox{ for all }U,V\in TM,
\end{equation}
or equivalently $
\cur (Y_1,Y_2)=0 \mbox{ for all } Y_1,Y_2\in X^\bot .$
\end{de}

Since $X$ is not parallel all the trace conditions which were true for a pp-wave,  fail to
hold for a pr-wave. 
For example, if we suppose (\ref{preq}) we get for the trace
$tr_{(3,5)(4,6)}(\cur \otimes \cur )(U,V,W,Z)=
h_p(\cur (U,V)X, \cur(W,Z)Z)$ which is not necessarily zero.
But one can prove an equivalence similarly to \ref{eq1} of
Lemma \ref{eqs}. (For the proof of the following statements see \cite{leistner05}.)

\blem
A Lorentzian manifold $(M,h)$ with recurrent vector field $X$ is a  
pr-wave if and only if 
$ \Lambda_{(1,2,3)} \left(\xi \otimes {\cal R}\right) =0$, where $\xi$ denotes the
1-form $h(X,.)$.
\elem

Also we get a similar description in terms of local coordinates as for pp-waves.
\blem 
A Lorentzian manifold $(M,h)$ of dimension $n+2>2$ is a pr-wave if and only if 
around any point $o\in M$ exist coordinates
 $(U,\varphi=(x, (y_i)_{i=1}^n, z)) $ in which the metric $h$ has the following
form,
\begin{equation}
\label{prform}h = 2\ dx dz   + f dz^2 +  \sum_{i =
1}^n  dy_i^2  \ \mbox{, with $ f\in C^\infty(M)$.}
\end{equation}

\elem

Regarding the vanishing of the screen holonomy  the following result can be obtained by the description of proposition \ref{ppeinfach1} and the definition of a pr-wave.
\bs
A  Lorentzian manifold  $(M,h)$ with recurrent vector field is a pr-wave if and only if the following equivalent conditions hold:
\bnum
\item[(1)]
 The screen holonomy of
$(M,h)$ is trivial (i.e. the screen bundle over $M$ is flat).
\item[(2)] $(M,h)$ has solvable holonomy contained in $\rr\ltimes \rrn$.
\enum
In addition, $(M,h)$  is a pp-wave if and only if its holonomy is Abelian, i.e. contained in $\rrn$.
\es

Finally, we see that Ricci-isotropy forces a pr-wave to be a pp-wave.

\bs
A pr-wave is a pp-wave if and only if it is Ricci-isotropic.
\es

For sake of completeness we also want to mention two subclasses of pp-waves. 
The first are the plane waves which are pp-waves with quasi-recurrent curvature, i.e. 
$\nabla \cur = \xi \otimes \tilde{\cur}$ where $\xi=h(X,.)$ and $\tilde{\cur}$ a
$(4,0)$-tensor.
For plane waves the function $f$ in the 
local form of the metric is of the form $f=\sumij a_{ij} y_i y_j$ where the $a_{ij}$ are functions of $z$.
A subclass of plane waves are the Lorentzian symmetric spaces with solvable transvection group, the so-called Cahen-Wallach spaces (see \cite{cahen-wallach70}, also \cite{bb-ike93}). Here the function $f$ satisfies
 $f=\sumij a_{ij} y_i y_j$ where the $a_{ij}$ are constants.

\section{Another generalisation of pp-waves}
\label{secrr}
Now we introduce another class of Lorentzian manifolds by relaxing also the curvature condition. 

\begin{de}\label{rrdef}
We say that a Lorentzian manifold $(M,h)$ with recurrent vector field has {\em light-like hypersurface curvature} if its curvature tensor ${\cal R}$ obeys
\begin{equation}
\label{rreqn}
\cur (U,V): \Xi^\bot \longrightarrow \Xi \mbox{ for all }U,V\in \Xi^\bot
\end{equation}
where $\Xi$ and $\Xi^\bot$ are the light-like distributions defined by the recurrent vector field.
\end{de}
Of course, (\ref{rreqn}) is equivalent to the fact the the $(4,0)$-curvature tensor 
vanishes on $\Xi^\bot \times\Xi^\bot \times \Xi^\bot \times \Xi^\bot$.

The chosen name can be explained by the following considerations. Since $(M,h)$ carries a recurrent vector field, the manifold is foliated into the flow of this vector field and the submanifolds defined by the integrable distribution $\Xi^\bot$. Hence, through any point $p\in M$ goes a one-dimensional isotropic submanifold $\cal X_p$ and a light-like hypersurface $\cal X_p^\bot$ with tangent bundles  $T\cal X_p=\Xi|_{\cal X_p}$ and $T\cal X_p^\bot=\Xi^\bot|_{\cal X_p^\bot}$ respectively, satisfying $\cal X_p\subset \cal X_p^\bot$. Since the distribution $\Xi^\bot$ is parallel, i.e. $\nabla_U:\Gamma(\Xi^{(\perp)})\rightarrow \Gamma(\Xi^{(\perp)})$
for every $U\in TM$, the Levi-Civita connection $\nabla$ of $(M,h)$ defines a connection on the hypersurface $\cal X_p^\bot$, denoted by $\ring{\nabla}: \Gamma ( T \cal X_p^\bot \otimes T\cal X_p^\bot )\rightarrow \Gamma( T\cal X_p^\bot )$. 
Then we get the following equivalences.
\bs
A Lorentzian manifold with recurrent vectorfield $X$ has light-like hypersurface curvature if and only if every light-like hypersurface $\cal X_p^\bot$, defined by $X$ and equipped with induced connection $\ring{\nabla}$, satisfies one of the following equivalent conditions:
\bnum
\item[(1)] $\ring{\cal R}(U,V)W$ is light-like, for $U,V,W\in T\cal X_p^\bot$ and $\ring{\cal R}$ the curvature of $\ring{\nabla}$.
\item[(2)] The holonomy of $\ring{\nabla}$ is solvable and contained in $\rr\ltimes \rrn$. 
\enum
If in addition $X$ is parallel, then the holonomy of $\ring{\nabla}$ is Abelian and  contained in $\rrn$.
\es
\bprf First we prove the equivalence of both conditions under the assumption that $(M,h)$ is a Lorentzian manifold with recurrent vector field $X$. 
The equivalence is based on  the  Ambrose-Singer holonomy theorem which says that $\mf{hol}_q(\cal X_p^\bot, \ring{\nabla})$ is generated the following endomorphisms  of $T_q\cal X_p^\bot=\Xi_q^\bot$,
\[ \ring{\cal P}^{-1}_\gamma \circ \ring{\cal R}(U,V) \circ \ring{\cal P}_\gamma \in \lagl(T_q\cal X_p^\bot),
\]
where $\ring{\cal P}_\gamma$ is the parallel displacement w.r.t. $\ring{\nabla}$ along a curve $\gamma$ in $\cal X_p^\bot$ starting at $q$, and $U,V\in \Xi^\bot_{\gamma(1)}$. For $\gamma$ the constant curve it becomes evident that (2) implies (1). Now bearing in mind that $\nabla=\ring{\nabla}$ on $\cal X_p$ which implies that $\ring{\cal P}_\gamma$ leaves $\Xi$ invariant, (1) implies that the holonomy algebra maps $T_q\cal X_p^\bot=\Xi_q^\bot$ onto $\Xi_q$ which means it is contained in
\[ \left\{\left.\left(\begin{array}{cc}a&v^t\\0&0\end{array}\right)\ \right|\ a\in \rr, v\in \rrn\right\}\subset \lagl(n+1)\]
with respect to a basis adapted to $\Xi_q\subset \Xi_q^\bot$.  In addition, when $X$ is parallel it is mapped  to zero by the holonomy algebra as $\ring{\cal R}(U,V)X=0$.
Finally it is evident that the condition (\ref{rreqn}) from the definition is equivalent to (1).
\eprf
In \cite{bezvitnaya05} the quantities assigned to the hypersurfaces $\cal X_p^\bot$ are used to describe the holonomy of a Lorentzian manifold further, in particular to decide to which type in the distinction following Berard-Bergery and Ikemakhen in \cite{bb-ike93} the holonomy algebra belongs. This approach makes use of a {\em screen distribution} which is  complementary and orthogonal to $\Xi$ in $\Xi^\bot$ (see also \cite{bejancu-duggal}). Such a screen distribution can always be chosen, but since it requires a choice we prefer to work with an analogon to the screen bundle introduced in section \ref{secscreen} which can be defined without making such a choice.

Let $\cal X_p^\bot$ be a light-like hypersurface through $p\in M$ defined by the recurrent vector field $X$. Then we define the {\em restricted screen bundle} over $\cal X_P^\bot$ as
\[\ring{\cal S}\ :=\ {\cal S}|_{\cal X_p^\bot}.\]
$\ring{\cal S}$ is equipped with a covariant derivative defined by $\ring{\nabla}$,
\[\nabla^{\ring{\cal S}}_U [V]\ :=\ \left[ \ring{\nabla}_U V\right],\ \text{ for } U\in T\cal X_p^\bot, V\in \Gamma (T\cal X_p^\bot).\]
Again, since $\Xi$ is $\ring{\nabla}$-invariant, this is well-defined. We obtain another equivalence in terms of the screen bundle.
\bs
A Lorentzian manifold with recurrent vectorfield $X$ has light-like hypersurface curvature if and only if over every light-like hypersurface $\cal X_p^\bot$ defined by $X$ the connection $\nabla^{\ring{\cal S}}$ on the restricted screen bundle $\ring{\cal S}$  is flat.\es

\bprf 
The curvature ${\cal R}^{\ring{\cal S}}$  of $\nabla^{\ring{\cal S}}$ can be written in terms of the curvature of $\ring{\nabla}$ as
\[{\cal R}^{\ring{\cal S}}(U,V)[W]\ =\ \left[\ring{\cal R}(U,V)W\right] ,\text{ for }U,V,W\in T\cal X_p^\bot.\]
Then the previous proposition gives the equivalence.
\eprf

For the case where the vector field $X$ is parallel we obtain the following equivalent trace condition.

\bs
A Brinkmann wave $(M,h)$ has light-like hypersurface curvature if and only if the curvature tensor $\cal R$ of $(M,h)$ obeys
$||\cur ||^2=0$ where 
$||\cur ||^2$ is the square of the norm of the curvature tensor, defined by
$||\cur ||^2:= tr_{(1,5)(2,6)(3,7)(4,8)}(\cur \otimes \cur )  $.
\es
\bprf
Again we use the basis $(X,E_1, \ldots ,E_n,Z)$ as in (\ref{bundlep}). Because $X$ is parallel every curvature term where
$X$ is plugged in vanishes and we get
\[||\cur ||^2=\sum_{i,j,k,l=1}^n \cur(E_i,E_j,E_k,E_l)^2.\]
But this expression vanishes if and only if $(M,h)$ satisfies (\ref{rreqn}).
\eprf

Now we want to focus on the description of a Lorentzian manifold with light-like hypersurface curvature in local coordinates.
\bs\label{rrformsatz}
A Lorentzian manifold $(M,h)$ of dimension $n+2>2$ has light-like hypersrface curvature if and only if 
around any point $o\in M$ exist coordinates
 $(U,\varphi=(x, (y_i)_{i^=1}^n, z)) $ in which the metric $h$ has the following
local shape
\begin{equation}
\label{rrform}h = 2\ dx dz   + f dz^2  
+ \left(\sumi u_i\ dy_i\right)dz+
 \sum_{i=1}^n  dy_i^2  
\end{equation}
with $\frac{\partial u_i }{\partial x}=0$ and
$ f\in C^\infty(M)$. If, in addition, $(M,h)$ is  a Brinkmann wave, then $f$ does not depend on $x$.
In the corresponding Schimming coordinates ($u_i=0$) the $g_{ij}$ are the coefficients of 
a $z$-dependent family of flat Riemannian metrics.
\es
\bprf
For Walker coordinates $(x,y_1, \ldots , y_n,z)$ the condition that  ${\cal R}$ vanishes on $\Xi^\bot$ gives that
$ {\cal R}(\ddi, \ddj, \ddk, \ddl)=0$. But this is the integrability condition for the existence of new coordinates
with $h(\ddi,\ddj)=\delta_{ij}$. This proves the first point.

If we now chose Schimming coordinates we still have the condition
$ {\cal R}(\ddi, \ddj, \ddk, \ddl)=0$. But for Schimming coordinates
$ {\cal R}(\ddi, \ddj, \ddk, \ddl)$ equals to ${\cal R}^g(\ddi, \ddj, \ddk, \ddl)$ 
where ${\cal R}^g$ denotes the curvature tensor of the Riemannian metrics defined by the coefficients
$g_{ij}$. Hence for each $z$ this has to be a flat Riemannian metric.   
\eprf
The description in these coordinates shows that the $\lason$-part of the curvature and the holonomy is generated by expressions of the form $\cal R(\ddz,\ddi)$ as we will see in the following.
We will illustrate this description in two different types of coordinates by some calculations.

First we calculate the curvature of a Lorentzian manifold with light-like hypersurface curvature in a point and  given  coordinates of the form 
(\ref{rrform}).
We can arrange these coordinates around the point $p$ in way that
 $\ddx, (\ddi)_{i=1}^n,\ddz$ is a basis of the form $(X,E_1, \ldots ,E_n,Z)$. Hence, if 
$\lag:=pr_{\lason}(\mf{hol}_p(M,h))$, then $\lag$ contains the following elements of $\lason$, for each
$U,V\in T_p M$:
\[\sumij h\left({\cal R}(U,V)\ddi,\ddj\right)E_{ij}\]
where $E_{ij}$ denotes the standard basis of $\lason$.
Now the only non-vanishing curvature terms of this form are 
\[h\left({\cal R}\left(\ddk,\ddz\right)\ddi,\ddj\right)=
{\cal R}\left(\ddk,\ddz,\ddi,\ddj\right)=\einhalb \ddk\left(\ddi (u_j) - \ddj( u_i)\right).\]
Hence, if one finds functions $(u_1, \ldots u_n)$ with 
$ \ddk\left(\ddi (u_j) - \ddj (u_i)\right)\not= 0$ one obtains a 
non-irreducible, non indecomposable Lorentzian manifold with light-like curved hypersurfaces, 
but with non-trivial
screen holonomy.

Now we calculate the curvature of such a manifold in Schimming coordinates, i.e. with $u_i=0$, i.e.
\[ h= 2dxdz +f dz^2  + \sumij g_{ij}dy_i dy_j.\]
Having light-like hypersurface curvature implies that $\cur(\ddi,\ddj,\ddk, \ddl)=0$, i.e. that $g_{ij}$ is a
$z$-dependent family of flat Riemannian metrics. 
Lets denote by $\Gamma^k_{ij}$ its $z$-dependent Christoffel symbols.
Then we get the following 
for the only non-vanishing curvature terms which are relevant for the $\lason$-projection of the holonomy:
\be
\cur (\ddi,\ddj,\frac{\partial}{\partial y_p},\ddz)
&
=
&
\einhalb \ddz \left( \ddj (g_{ip}) -\ddi (g_{jp}) \right)\\
&&{}+ \einhalb
\sumk
 \left(
\Gamma^k_{ip} \ddz (g_{jk}) -\Gamma^k_{jp} \ddz (g_{ik})\right)
\ee
In order to construct a non-irreducible, indecomposable Lorentzian manifold
 with lightlike curved hypersurfaces and non-trivial screen holonomy one has to find a family
of flat Riemannian metrics with Christoffel symbols such that the above expression is non-zero.

Now we prove to further properties of the coordinates.

\bs
A Lorentzian manifold with light-like hypersurface curvature is a pr-wave, i.e. has trivial screen holonomy, if 
there exist local coordinates of the form (\ref{rrform}) such that the $z$-dependent family of one-forms 
$\phi:=\sumk u_k d y_k$ on $\rrn$ is closed for any $z$.
\es
\bprf
Since $\phi_z$ are closed
--- considered as a family of differential forms on $\rrn$ --- 
they are a differential of a
function $\varphi$ which does not depend on the $x$ coordinate. More exactly: If $\phi_z=
\sumk u_k dy_k$ with $\ddx(u_k)=0$ and
\[0=d\phi_z= \sums_{l=1}^n du_l\wedge dy_l =\sums_{k,l=1}^n\ddk (u_l) dy_k\wedge dy_l\]
then exists a $\beta \in C^\infty(M)$ with $\ddx(\beta)=0$ and
$u_k= \ddk( \beta)$.

Now we consider the following coordinates
\beq\label{phitrafo}
\tilde{x}=x+\beta\ ,\;\tilde{y}_i=y_i\ ,\;\tilde{z}=z.
\eeq
These satisfy $\tilde{u}_i=0$, $\tilde{g}_{ij}=\delta_{ij}$ and $\tilde{f}=f-2 \ddz (\beta )$, and are therefore coordinates of a pr-wave.
\eprf

Regarding the Ricci curvature we can prove the condition for the Ricci isotropy in terms of the form $\phi$.

\bs
An Brinkmann wave with light-like hypersurface curvature is Ricci isotropic if and only if there are coordinates for which the family of one-forms
$\phi:=\sumk u_k d y_k$ on $\rrn$ satisfies the equation
\beq\label{lame}
d^*d\ \phi_z=0.
\eeq
for all $z$.
\es
\bprf We consider $\phi_z$ as a family of $1$-forms on  $\rrn$.
Fixing coordinates of the form (\ref{rrform}) we get that the basis
\be
X &:=& \deri{\partial x}\\
Z&:=&\deri{\partial z}- \frac{f}{2} \ddx\\
E_k&:=& \ddi - u_i\ddx.
\ee
is of the form (\ref{bundlep}). In these coordinates and this basis we obtain as $\ddx$ is parallel:
\be
Ric(\ddz,\ddi)
&=&\overbrace{{\cal R}((\ddz,X,Z, \ddi)}^{=0}+ \overbrace{{\cal R}((\ddz,X,Z, \ddi)}^{=0}+
\sumk {\cal R}\left(\ddz,\ddk,\ddk, \ddi\right)\\
&=&-\frac{1}{2}\sumk\left[ \ddk\Big(\ddk\big(\phi_z(\ddi)\big)\Big)- \ddk\Big(\ddi\big(\phi_z(\ddk)\big)\Big)\right]\\
&=&d^*d\phi_z(\ddi).
\ee
Here  $d^*$ is the co-differential with respect to the flat Riemannian metric
$g\equiv\delta_{ij}$. But a Brinkmann wave is Ricci-isotropic if and only if
$Ric(Y,.)=0$ for every $Y\in \Xi^\bot$ (see for example \cite{leistner05}) which gives the statement.
\eprf
The Ricci-isotropy is an important property because it is a necessary condition of the existence of parallel spinors on $(M,h)$.

\section{Further remarks on holonomy and examples}
\label{sechol}
We want to start the concluding remarks about the holonomy of Lorentzian manifolds with light-like hypersurface curvature with an example.
\begin{bsp}\label{ike}
There are examples which show that Lorentzian manifolds with light-like hypersurface curvature can have non-trivial screen holonomy, in particular having irreducible screen holonomy 
$\laso(3)\subset \laso(5)$ given by the Riemannian symmetric pair. The first example of such a manifold was 
given in \cite{ike96}, although with another purpose.  One considers  the following one-form
$\phi=\sum_{k=1}^5 u_k d y_k$ on $\rr^5$ with
\be
u_1&=&-y_3^2-4y_4^2-y_5^2,\\
u_3&=&-2\sqrt{3}\ y_2y_3-2y_4y_5,\\
 u_5&=&2\sqrt{3}\ y_2y_5+2y_3y_4,\\
 u_2&=&u_4\ =\ 0.
  \ee
 Now one defines the Lorentzian metric on $\rr^7$ by
\[ h:= 2 dx dz + f dz^2 + \phi dz +
\sum_{k=1}^5 dy_k^2\]
where $f$ is a function on $\rr^7$ with $\frac{\partial f}{\partial y_i}\not=0$.
 The holonomy of this manifold equals
to $(\rr\+ \laso(3,\rr))\ltimes \rr^5$ or if $f$ does not depend on
$x$ equal to $\laso(3,\rr)\ltimes \rr^5$ where
$\laso(3,\rr)\subset \laso(5,\rr)$ is the irreducible representation defined by the Riemannian symmetric pair:
the Lie algebra $\mf{sl}(3,\rr)$ can be decomposed into vector spaces
$\mf{sl}(3,\rr)=\laso(3,\rr)\+ sym_0(3,\rr)$,
where $ sym_0(3,\rr)$ denote the trace free symmetric matrices. This is a $5$--dimensional
vector space, invariant and irreducible under the adjoint action of $\laso(3,\rr)$.
This representation is equal to the holonomy representation of the Riemannian symmetric space 
$Sl(3,\rr)/SO(3,\rr)$.

Another example of this type having the same holonomy was constructed in \cite{thesis} by setting
$
u_1= -4 y_1y_2$,
$u_2= 4 y_1y_2$,
$u_3= - y_1y_4 - y_2y_4+ y_1y_3 -y_2y_3 + \sqrt{3}(y_4y_5 -y_3y_5)$,
$u_4=  y_1y_4 - y_2y_4+ y_1y_3 +y_2y_3 + \sqrt{3}(y_4y_5 +y_3y_5)$ and 
$u_5=0$.
Recently in \cite{galaev05} another such example was  constructed by defining
$
u_1=-\frac{2}{3}((y_3)^2+4(y_4)^2+(y_5)^2)$,
$ u_2=\frac{2\sqrt{3}}{3}((y_3)^2-(y_5)^2)$,
$u_3=\frac{2}{3}(y_1y_3-\sqrt{3}y_2y_3-3y_4y_5-(y_5)^2)$,
$ u_4=\frac{8}{3}y_1y_4$ and $
u_5=\frac{2}{3}(y_1y_5+\sqrt{3}y_2y_5+3y_3y_4+y_3y_5)$.
These examples also have $\laso(3)\subset \laso(5)$ as screen holonomy.
We do not know whether the three examples are locally isometric.

On the other hand one can construct a manifold with the same holonomy but with different geometric properties, i.e.
which does not have light-like hypersurface curvature, by the following construction. Let $g$ be the Riemannian metric 
on $Sl(3,\rr)/SO(3,\rr)$ and consider the Lorentzian manifold
\[ \left( M:=\rr^2\times Sl(3,\rr)/SO(3,\rr), h:=2dxdz +f dz^2 + g.\right)\]
If $f$ is sufficient general this manifold is indecomposable and has holonomy algebra
$\laso (3)\ltimes \rr^5$ or $(\rr \+ \laso (3)) \ltimes \rr^5$. But, its curvature restricted to $\left(
\ddx\right)^\bot$ does not vanish because $\cur\left(\ddi,\ddj,\ddk,\ddl\right)$ equals to the curvature
of $Sl(3,\rr)/SO(3,\rr)$.  
\end{bsp}
\vspace{.5cm}
This example as well as the curvature calculations in local coordinates show that the $\lason$-part of the 
holonomy of a Lorentzian manifold with light-like hypersurface curvature is not necessarily trivial. Due to a recent result in \cite{galaev05} one can even show that any possible screen holonomy, i.e. any Riemannian holonomy group can occur as screen holonomy of a Lorentzian manifold with light-like hypersurface curvature.  We will now indicate why this is the case.

The classification of possible screen holonomies was based on the notion of {\em weak curvature endomorphisms} and {\em weak-Berger algebras} which were introduced in \cite{leistner02}. Weak curvature endomorphism are defined for a Lie algebras $\lag\subset \lason$ by a Bianchi-identity:
\[
\cal B(\lag) := \left\{Q\in Hom(\rrn,\lag)\mid \la Q(x)y,z\ra + \la Q(y)z,x\ra+\la Q(z)x,y\ra\right\}
\]
$\cal B(\lag)$ is the kernel of the homomorphisms $\lambda:Hom(\rrn,\lag)\rightarrow \Lambda^3(\rrn)^*$ which is the combination of skew symmetrisation and dualisation by means of $\la.,.\ra$. $\lag$ is called {\em weak-Berger algebra} if and only if
\[\lag\ =\ \mathsf{span} \{Q(x)\mid Q\in \cal B(\lag), x\in \rrn\}.\]
In \cite{leistner03} and \cite{leistner03b} we showed that any weak-Berger algebra is a Riemannian holonomy algebra. 
On the other hand, any Riemannian holonomy algebra can be realised as screen holonomy of a Lorentzian manifold with recurrent or parallel light-like vector field by the following construction. Let $(N,g)$ be a Riemannian manifold and $f\in C^\infty (N\times \rr) $ a smooth function which is sufficiently generic. Then $M:=\rr^2\times N$ with the metric
\[h\ :=\  2dxdz + f dz^2 + g\]
is a Lorentzian manifold with recurrent vector field and the screen holonomy of $(M,h)$ is equal to the Riemannian holonomy $Hol_p(N,g)$ (see \cite{leistner01}). But it was an open question if for any of the four types of holonomy  in \cite{bb-ike93} any Riemannian holonomy can be realised as screen holonomy. In \cite{galaev05} it was shown that  this is possible. We will now describe briefly parts of this method which we will need to construct further examples. This construction uses the fact that the screen holonomy $\lag$ is a weak-Berger algebra. For details of the following see \cite{galaev05}.

First, for a weak-Berger algebra $\lag\subset \lason$ one fixes weak curvature endomorphisms $Q_A\in \cal B(\lag)$ for $A=1,\ldots, N$ and a basis $e_1, \ldots , e_n$ of $\rrn$, orthonormal w.r.t. $\la.,.\ra$.  Now one defines the following polynomials on $\rr^{n+1}$,
\beq\label{poly}
u_i(y_1, \ldots , y_n,z)\ :=\ \sum_{A=1}^N \sumkl
\frac{(A-1)!}{3}\Big\la Q_A(e_k)e_l +  Q_A(e_l)e_k,e_i\Big\ra y_k y_l z^A,\eeq
and the following Lorentzian metric on $\rr^{n+2}$
\beq\label{polymetric} h\ =\ 2 dx dz + f dz^2 +2 \sumi u_idy_idz +\sumk dy_k^2,\eeq
where $f$ is a function on $\rr^{n+2}$.
This metric is analytic, hence its holonomy is generated by the derivations of the curvature tensor. But the metric is constructed in a way such that the only non-vanishing $\lason$-parts of the curvature and its derivatives satisfy
\beq\label{polycurv}pr_{\lason}\Bigg[ \big(
\underbrace{\nabla_{\ddz}\ldots\nabla_{\ddz}}_{(A-1)-\text{times}}\cal R\big)\left(\ddi,\ddz\right)\Bigg]\ =\ Q_A(e_i),\eeq
 for $ A=1, \ldots , N$ and $i=1,\ldots ,n$. If one now starts this construction with $Q_1, \ldots, Q_N$ which span $\cal B(\lag)$, e.g. a basis of $\cal B(\lag)$, then the derivatives of the curvature will span 
$\lag$. Hence the weak-Berger algebra $\lag$ we started with is the screen holonomy of $(\rr^{n+2},h)$. 
But, more importantly, it is proven that, for any of the four types of indecomposable, non-irreducible Lorentzian holonomy in \cite{bb-ike93}  the function $f$ can be chosen in a way that the holonomy of $h$ belongs to this type. 

For our purposes it is important that the constructed metric $h$ admits light-like hypersurface curvature due to the description in coordinates in  proposition \ref{rrformsatz}. We obtain the following result.
\bs
For any of the four types of indecomposable, non-irreducible Lorentzian holonomy and any Riemannian holonomy algebra $\lag$ there is a Lorentzian manifold $(M,h)$ with light-like hypersurface curvature such that the holonomy of $(M,h)$ is of the given type and its screen holonomy is equal to $\lag$.
\es
This result is most remarkably as the curvature condition on a manifold with light-like hypersurface curvature are very strong and only a slight generalisation of the curvature conditions posed on a pp-wave.

\bigskip

The method described above gives a construction principle for Lorentzian manifolds with light-like hypersurface curvature under the assumption that the weak curvature endomorphism are known. But since every weak Berger algebra $\lag$ is a Riemannian holonomy algebra and thus a Berger algebra, i.e. 
\[\lag = \mathsf{span}\{R(x,y)\mid R\in \cal K(\lag), x,y\in \rrn\}\]
where $\cal K(\lag)$ are the following curvature endomorphisms
\[\cal K(\lag)=\{R\in Hom(\Lambda^2\rrn, \lag)\mid R(x,y)z+R(y,z)x+R(z,x)y=0\},\]
sometimes it is sufficient to know the space $\cal K(\lag)$. We will illustrate this in the following construction which generalises  Example \ref{ike}.
First we we note that both spaces of curvature endomorphisms, $\cal B(\lag)$ and $\cal K(\lag)$ are $\lag$-modules and their relation is as follows.

\blem\label{lemma1} Let $\lag\subset \lason$. Then
the vector space ${\cal R}(\lag)$  spanned by $ \{ R(x,.)\in 
{\cal B}(\lag) \ |\ R\in \kg,\ x \in \rrn\}$ is
a $\lag$-submodule of $\cal B(\lag)$.
\elem
\bprf
Because of the defining Bianchi-identities $\cal R(\lag)\subset \cal B(\lag)$ is ensured. 
For $R\in \cal K(\lag)$ it is 
\[\left(A \cdot R(x,.)\right) (y) \ =\    [A, R(x,y)] - R(x,Ay)\ =\ (A\cdot R)(x,y)+ R(Ax,y)\  \in\  \cal R(\lag),\]
which shows that $\cal R(\lag)$ is also a submodule.
\eprf
This lemma shows that apriori any Berger algebra is a weak-Berger algebra whereas the other implication requires a proof based on representation theory (see \cite{leistner02}, \cite{leistner03} and \cite{leistner03b}).
Nevertheless we can apply the lemma in order to construct examples of Lorentzian manifolds with light-like hypersurface curvature and the screen holonomy of a Riemannian symmetric space $G/K$.

Suppose $G/K$ is a  Riemannian symmetric space of semisimple type and of dimension $n$. In particular, the Lie algebras satisfy $\lag=\mf{k}\+\mf{m}$ with $\lak$ a subalgebra acting irreducible on $\lam$ and $[\lam,\lam]\subset \lak$. The metric on $G/K$ corresponds to an invariant inner product  $\la.,.\ra$ which is a multiple of the Killing form  $B$ of $\lag$. The holonomy group of $G/K$ is $K$ acting by the adjoint representation on $\lam\simeq T_{[e]} \left(G/K\right)$.
Suppose $X_1, \ldots X_n$ is a basis of $\lam$ which is orthogonal with respect to  $B$. Using these ingredients we define the following polynomials on $ \rr^{n+1}$:
\be
\lefteqn{u^{(G,K)}_i(y_1, \ldots , y_n,z)\ :=}\\
&&\sum_{j,k,l=1}^n 
\frac{(j-1)!}{3}\Big( 
B\big( [X_j,X_k],[X_l,X_i]\big)+
B\big( [X_j,X_l],[X_k,X_i]\big)
\Big)y_k y_l z^j,\ee
where $[.,.]$ is the commutator in $\lag$ and $B$ the Killing form.
Again we define a Lorentzian metric on $\rr^{n+2}$ by
\[ h^{(G,K)}\ =\ 2 dx dz + f dz^2 +2 \sumi u^{(G,K)}_idy_idz +\sumk dy_k^2,\]
for $f$ being a smooth function on $\rr^{n+2}$. In this situation it holds the following proposition.
\bs
Let $G/K$ be an irreducible Riemannian symmetric space of dimension $n$. Then the Lorentzian metric $h^{(G,K)}$ on $\rr^{n+2}$ has light-like hypersurface curvature and its screen holonomy is equal to the holonomy of the Riemannian symmetric space $G/K$, i.e. is equal to $K$.\label{symprop}
\es

\bprf The proof relies on the method described above.
Since $G/K$ is a Riemannian symmetric space, the curvature endomorphisms of $\lak$ satisfy $\cal K(\lak)=\rr\cdot [.,.]$, where $[.,.]$ is the commutator of $\lag$. Since $\lak$ is the holonomy algebra of this space we get $\lak=\mathsf{span}\{ [X,Y]\mid X,Y\in \lam\}$. Hence for a basis $X_1,\ldots , X_n$  of $\lam$, the $Q_j:=[X_j.,]$ span the submodule $\cal R(\lak)$ in $\cal B(\lak)$ by Lemma \ref{lemma1} and generate the whole Lie algebra $\lak$. In this situation, if the basis $X_i$ is assumed to be orthogonal w.r.t. the Killing form $B$, we obtain for the terms in (\ref{poly})
\be
B\big( Q_j(X_k)X_l +Q_j(X_l)X_k\ ,\ X_i\big) &= &
B\big( \left[X_j,X_k\right]X_l +\left[X_j,X_l\right]X_k\  ,\ X_i\big)\\
&=& 
B\big(  \left[ \left[X_j,X_k\right],X_l\right] +\left[\left[X_j,X_l\right],X_k\right]\ ,\ X_i\big)\\
&=&
B\big(   \left[X_j,X_k\right],\left[X_l,X_i\right]\big)+B\big(\left[X_j,X_l\right],\left[X_k,X_i\right]\big).
\ee
Hence, the curvature of $h^{(G,K)}$ satifies (\ref{polycurv}) which implies that the holonomy of $h^{(G,H)}$ is equal to $K$.
\eprf
Again, as in Example \ref{ike}, a Lorentzian manifold with the same screen holonomy can be obtained by the metric $h=2dxdz+fdz^2=g$ where $g$ is the Riemannian metric of $G/K$. But this manifold does not have light-like hypersurface curvature and is therefore not isometric to $h^{(G,K)}$.

In principle, the method of \cite{galaev05} works for any Riemannian holonomy algebra, also non-symmetric ones, if one is able to calculate $\cal B(\lag)$. As in Proposition \ref{symprop} one could also try to use the sub-module $\cal R(\lag)$, but for non-symmetric Riemannian holonomy groups $\cal K(\lag)$ can be very big and thus the calculations complicated. Another way is to use other, easier submodules of $\cal B(\lag)$. This methods works if $\lag$ is simple, since any sub-module of $\cal B(\lag)$ generates a non-trivial ideal in $\lag$ which has to be $\lag$ in this case. For example in the case of $\lag_2\subset \laso(V)$ with $V=\rr^7$ the $\lag_2$-module $Hom(V, \lag_2)$ which contains $\cal B(\lag_2)$ splits into the direct sum of 
$V_{[1,1]}$, $\odot^2_0V^* $ and $V$ where $V_{[1,1]}$ is the $64$-dimensional $\lag_2$-module of highest weight $(1,1)$ and $\odot^2_0V^* $ the $27$-dimensional module of  highest weight $(2,0)$. Since $\cal B(\lag_2)$ is the kernel of the skew-symmetrisation
\barr{rcrcl}
\lambda&:&Hom(V,\lag_2)&\rightarrow& \Lambda^3V^*\\
&& // \hspace{.8cm}&&\hspace{.3cm} \setminus\setminus \\[.2cm]
&&V_{[1,1]}\+\odot^2_0V^*\+V&&\odot^2_0V^*\+V\+\ccc
\earr
a dimension analysis shows that $\cal B(\lag)$ must contain $V_{[1,1]}$. Thus, by chosing a basis of $V_{[1,1]}$ a metric of the form (\ref{polymetric}) with coefficients as in (\ref{poly}) can be defined and one obtains a Lorentzian manifold with light-like hypersurface curvature and screen holonomy $G_2$.

\renewcommand{\baselinestretch}{1}


\begin{thebibliography}{Lei03b}

\bibitem[Bau02]{baumsurvey}
Helga Baum.
\newblock {Conformal Killing spinors and special geometric structures in
  Lorentzian geometry - a survey}.
\newblock In {\em Proceedings of the Workshop on Special Geometric Structures
  in String Theory, Bonn, September 2001}, 2002.
\newblock Proceedings archive of the EMS Electronic Library of Mathematics,
  {\tt www.univie.ac.at/EMIS/proceedings/}.

\bibitem[BBI93]{bb-ike93}
Lionel B{\'e}rard-Bergery and Aziz Ikemakhen.
\newblock On the holonomy of {L}orentzian manifolds.
\newblock In {\em Differential Geometry: Geometry in Mathematical Physics and
  Related Topics (Los Angeles, CA, 1990)}, volume~54 of {\em Proc. Sympos. Pure
  Math.}, pages 27--40. Amer. Math. Soc., Providence, RI, 1993.

\bibitem[BD96]{bejancu-duggal}
Aurel Bejancu and Krishan~L. Duggal.
\newblock {\em Lightlike Submanifolds of Semi-{R}iemannian Manifolds and
  Applications}, volume 364 of {\em Mathematics and Its Applications}.
\newblock Kluwer Academic Press, 1996.

\bibitem[Bez05]{bezvitnaya05}
Natalia Bezvitnaya.
\newblock {Lightlike foliations on Lorentzian manifolds with weakly irreducible
  holonomy algebra}, 2005.
\newblock arXiv:math.DG/0506101.

\bibitem[Bou00]{boubel00}
Charles Boubel.
\newblock {\em Sur l'holonomie des vari\'{e}t\'{e}s pseudo-riemanniennes}.
\newblock PhD thesis, Universit\'{e} Henri Poincar\'{e}, Nancy, 2000.

\bibitem[Bri25]{brinkmann25}
H.~W. Brinkmann.
\newblock {E}instein spaces which are mapped conformally on each other.
\newblock {\em Math. Ann.}, 94:119--145, 1925.

\bibitem[CW70]{cahen-wallach70}
M.~Cahen and N.~Wallach.
\newblock {L}orentzian symmetric spaces.
\newblock {\em Bull. Amer. Math. Soc.}, 79:585--591, 1970.

\bibitem[Gal05]{galaev05}
Anton~S. Galaev.
\newblock {Metrics that realize all types of Lorentzian holonomy algebras},
  2005.
\newblock arXiv:math.DG/0502575.

\bibitem[Ike96]{ike96}
Aziz Ikemakhen.
\newblock Examples of indecomposable non-irreducible {L}orentzian manifolds.
\newblock {\em Ann. Sci. Math. Qu\'{e}bec}, 20(1):53--66, 1996.

\bibitem[Kat99]{kathhabil}
Ines Kath.
\newblock {\em Killing Spinors on Pseudo-Riemannian Manifolds}.
\newblock 1999.
\newblock Habilitationsschrift, Humboldt-Universit\"at Berlin.

\bibitem[Lei02a]{leistner02}
Thomas Leistner.
\newblock Berger algebras, weak-{B}erger algebras and {L}orentzian holonomy,
  2002.
\newblock SFB 288-Preprint no. 567, {\tt ftp://
  ftp-sfb288.math.tu-berlin.de/pub/Preprints/preprint567.ps.gz}.

\bibitem[Lei02b]{leistner01}
Thomas Leistner.
\newblock {L}orentzian manifolds with special holonomy and parallel spinors.
\newblock In {\em The proceedings of the 21st winter school "Geometry and
  Physics", Srni, January 13-20, 2001}, number~69 in Supplemento ai Rendiconti
  del Circolo Matematico di Palermo. Serie II, pages 131--159, 2002.

\bibitem[Lei03a]{leistner03}
Thomas Leistner.
\newblock {Towards a classification of Lorentzian holonomy groups}, 2003.
\newblock arXiv:math.DG/0305139.

\bibitem[Lei03b]{leistner03b}
Thomas Leistner.
\newblock {Towards a classification of Lorentzian holonomy groups. Part II:
  Semisimple, non-simple weak Berger algebras}, 2003.
\newblock arXiv: math.DG/0309274.

\bibitem[Lei04]{thesis}
Thomas Leistner.
\newblock {\em Holonomy and Parallel Spinors in Lorentzian Geometry}.
\newblock Logos Verlag, 2004.
\newblock Dissertation, Mathematisches Institut der Humboldt-Universit\"{a}t
  Berlin, 2003.

\bibitem[Lei05]{leistner05}
Thomas Leistner.
\newblock {Conformal holonomy of C-spaces, Ricci-flat, and Lorentzian
  manifolds}, 2005.
\newblock arXiv:math.DG/0501239, to appear in {\em Differential Geometry and
  its Applications}.

\bibitem[Sch74]{schimming74}
Rainer Schimming.
\newblock {R}iemannsche {R}\"{a}ume mit ebenfrontiger und mit ebener
  {S}ymmetrie.
\newblock {\em Mathematische Nachrichten}, 59:128--162, 1974.

\bibitem[Wal49]{walker49}
A.~G. Walker.
\newblock On parallel fields of partially null vector spaces.
\newblock {\em Quart. Journ. of Mathematics}, 20:135--145, September 1949.

\end{thebibliography}

{\sc School of Mathematical Sciences,}\\
{\sc University of Adelaide,}\\
{\sc SA 5005,}
{\sc Australia.}\\
Email: {\tt tleistne@maths.adelaide.edu.au}

\end{document}